\documentclass[11pt,a4paper,reqno]{amsart}

\usepackage[utf8]{inputenc}

\usepackage[top=2.7cm,bottom=2.7cm,left=2.5cm,right=2.5cm]{geometry}

\usepackage{amsfonts,amsmath,amssymb,amsthm,dsfont,amsxtra,enumerate,multicol,mathtools,tkz-berge,hyperref,graphicx}
\usepackage{listings}
\usepackage{xcolor}

\definecolor{codegreen}{rgb}{0,0.6,0}
\definecolor{codegray}{rgb}{0.5,0.5,0.5}
\definecolor{codepurple}{rgb}{0.58,0,0.82}
\definecolor{backcolour}{rgb}{0.95,0.96,0.98}

\lstdefinestyle{mystyle}{
    backgroundcolor=\color{backcolour},   
    commentstyle=\color{codegreen},
    keywordstyle=\color{magenta},
    numberstyle=\tiny\color{codegray},
    stringstyle=\color{codepurple},
    basicstyle=\ttfamily\footnotesize,
    breakatwhitespace=false,         
    breaklines=true,                 
    captionpos=b,                    
    keepspaces=true,                 
    numbers=left,                    
    numbersep=5pt,                  
    showspaces=false,                
    showstringspaces=false,
    showtabs=false,                  
    tabsize=2
}

\newtheorem{theorem}{Theorem}[section]
\newtheorem{proposition}[theorem]{Proposition}
\newtheorem{lemma}[theorem]{Lemma}
\newtheorem{corollary}[theorem]{Corollary}

\theoremstyle{definition}
\newtheorem{definition}[theorem]{Definition}
\newtheorem{remark}[theorem]{Remark}
\newtheorem{example}[theorem]{Example}
\newtheorem{notation}[theorem]{Notation}
\newtheorem{question}[theorem]{Question}

\newtheorem{alphconj}{Conjecture}

\author[Marie Amalore Nambi]{Marie Amalore Nambi}
\address{Sabanci University, Faculty of Engineering and Natural Sciences, Orta Mahalle, Tuzla, 34956, Istanbul, Turkey}
\email{amalore.p@gmail.com, amalore.pushparaj@sabanciuniv.edu}

\author[Neeraj Kumar]{Neeraj Kumar}
\address{Department of Mathematics, Indian Institute of Technology Hyderabad, Kandi, Sangareddy - 502285, INDIA}
\email{neeraj@math.iith.ac.in}

\title{On Conjecture of Binomial Edge Ideals of Linear Type}

\date{November 2023}

\begin{document}
\subjclass[2020]{{Primary 05E40, 13F65, 13F70}; Secondary {13A30}} 
%13A30 Associated graded rings of ideals (Rees ring, form ring), analytic spread and related topics
%13F65 Commutative rings defined by binomial ideals, toric rings, etc
%13F70 Other commutative rings defined by combinatorial properties
%05E40 Combinatorial aspects of commutative algebra
\keywords{Rees algebra, linear type, $d$-sequence, $p$-sequence, binomial edge ideal, tree.}

\date{}

\begin{abstract}
An ideal $I$ of a commutative ring $R$ is said to be of {\emph{linear type}} when its Rees algebra and symmetric algebra exhibit isomorphism. In this paper, we investigate the conjecture put forth by Jayanthan, Kumar, and Sarkar in \cite{JAR2021} that if $G$ is a tree or a unicyclic graph, then the binomial edge ideal of $G$ is of linear type. Our investigation validates this conjecture for trees. However, our study reveals that not all unicyclic graphs adhere to this conjecture. 
\end{abstract}

\maketitle

\section*{Introduction}

 Consider a graph $G$ on $[n]$ vertices. Let $S=\mathbb{K}[x_1,\ldots,x_n,y_1,\ldots,y_n]$ be a polynomial over a field $\mathbb{K}$. Herzog et al. \cite{HH} and independently Ohtani \cite{O2011} introduced the ideal $$J_G=(f_{ij} = x_iy_j-x_jy_i \mid i<j, \{i,j\} \in E(G)) \subset S$$ is called the \textit{binomial edge ideal} of the graph $G$. We call the binomial $f_{ij}$ associated to an edge $\{i,j\} \in E(G)$ as an \textit{edge binomial}. In \cite{HH}, Herzog et al. observed that the binomial edge ideal appears naturally in the study of conditional independence ideals, providing a valuable tool to investigate robustness theory in the context of algebraic statistics. In recent years, researchers have investigated various algebraic properties, including radicality, primality, Gr\"{o}bner bases, primary decomposition, among others, of binomial edge ideals in relation to their underlying combinatorial properties; see for partial list \cite{VRT2021, HH, MS, O2011, SK} and the references therein. 
 
\vspace{2mm}
 
 Let $A$ be a standard graded polynomial ring over a field $\mathbb{K}$ and $I \subset A$ be a homogenous ideal. The Rees algebra of an ideal $I \subset A$, is defined as $\mathcal{R}(I)= \{ \sum_{i=0}^n a_it^i \mid n \in \mathbb{N}$  and  $a_i \in I^i \} =\bigoplus_{i \geq 0}I^it^i \subset A[t]$, where $I^0=A$. $\mathcal{R}(I)$ is a classical algebraic structure that encodes numerous asymptotic properties of that ideal. 
 %Indeed, the Rees algebra provides insight into how the powers of an ideal sit inside another Noetherian ring.
 The symmetric algebra of an ideal $I \subset A$ is defined as $Sym(I)= T(I)/L$, where $T(I)$ is the tensor algebra of $I$ and $L$ is the ideal of $T(I)$ generated by elements of the form $x\otimes_A y- y \otimes_A x$, where $x,y \in I$. 
 Moreover, there exists a natural epimorphism from Sym$(I)$ to $\mathcal{R}(I)$, where Sym$(I)$ is the symmetric algebra of $I$. For an ideal $I\subset A$, if Sym$(I) \cong \mathcal{R}(I)$, then $I$ is referred to as an ideal of linear type. The defining equations of the symmetric algebra of an ideal can always be computed from a presentation matrix of the ideal (cf. \cite[Section 5.5]{HS2006}).

\vspace{2mm}
 
 Villareal provided a characterization of the ideal of linear type for edge ideals, namely, edge ideals are of linear type if and only if the graph is a tree or an odd unicyclic graph (cf. \cite{V1990}). This raises the question of whether a similar characterization might apply to other combinatorial objects. In \cite{JAR2021}, Jayanthan et al. initiated such a study for binomial edge ideals. In their work, the authors characterized almost complete intersection ideals and proposed the following conjecture:

 \begin{alphconj} \cite[Conjecture 4.17. (a)]{JAR2021} \label{conjlineartype}
    If $G$ is a tree or a unicyclic graph, then $J_G$ is of linear type.
\end{alphconj}

Huneke introduced the notion of $d$-sequence in \cite{H1980, H1982}. Huneke, and independently Valla, proved that if the ideal is generated by unconditioned $d$-sequence, then the ideal is of linear type in \cite{H1980, V1980}. Costa \cite{C1985} established that $d$-sequence condition on $I$ is a sufficient condition for the ideal $I$ to be of linear type. Later, Raghavan \cite{R1991} provided a simpler proof that $d$-sequence ideals are of linear type. To the best of our knowledge, study of the Rees algebra of binomial edge ideal of linear type are the following. In \cite{JAR2021}, the authors characterised almost complete intersection binomial edge ideals, a subclass of ideal of linear type. In \cite{AN, ANcycle}, the authors provided a partial answer for Conjecture \ref{conjlineartype}. Particularly, the authors characterized $d$-sequence binomial edge ideals of trees and unicyclic graphs at the level of a sequence of edge binomials. In \cite{A22}, Kumar characterized binomial edge ideals of linear type for closed graphs.

\vspace{2mm}

In this paper, our main aim is to settle Conjecture \ref{conjlineartype} for trees; see Theorem \ref{thm.conj}. To achieve this, we introduce a $p$-sequence in Section \ref{sec.p-seq}, as defined in Definition \ref{p-sequence}. 
%In (\cite[Corollary 6.2]{B2018}), Bolognini et al. proved that if $G$ is a bipartite graph, then the binomial edge ideal $J_G$ coincides with the parity binomial edge ideal $\mathcal{I}_G$ (cf. \cite{TCT}) and Lovász–Saks–Schrijver ideal $L_G(2)$ (cf. \cite{HMSW}), use which we establish the linear type properties of these ideals for trees. 
In Section \ref{sec.p-seqbei}, we proved that trees and a class of unicyclic graphs form a $p$-sequence, see Theorems \ref{thm.pseqtree} and \ref{thm.pseqcycle}. Finally, in Section \ref{sec.counterex}, we present a counterexample illustrating that the conjecture does not hold for unicyclic graphs in general; refer to Proposition \ref{prop.counter}.

\vspace{2mm}

%Secondly, we study normally torsionfree property of binomial edge ideals of trees. An ideal is said to be normally torsionfree if all its powers has same associated primes (cf. \cite[Definition 5.2]{SVV}). 

%Secondly, when $G$ is a tree, we study symbolic and ordinary powers of $J_G$, and the Cohen-Macaulay property of $\mathcal{R}(J_G)$. The study of symbolic powers is one of central topics in commutative algebra. In general symbolic powers and ordinary powers of an ideal do not coincide. There are some cases where it coincides. For the case of binomial edge ideals, Ohtani obtain the equality for complete multipartite graphs and Ene and Herzog obtain the equality for closed graphs (cf. \cite{EH, O2013}). In \cite{VRT2021}, the authors showed that $\mathcal{R}(J_G)$ is Cohen-Macaulay when $G$ is a closed graph. In this paper, we obtain the equality for trees by showing initial ideal of binomial edge ideals of trees are normally torsionfree, see Theorem \ref{thm.CM}. An ideal is said to be \textit{normally torsionfree} if all its powers has same associated primes (cf. \cite[Definition 5.2]{SVV}). Furthermore, we shown that the Rees algebra binomial edge ideals of trees are Cohen-Macaulay.

\vspace{2mm}

%This paper is organized as follows. The following section recalls definitions, notations, and known results required in the rest of the paper. In Section \ref{sec.p-seq}, we introduce and study $p$-sequence ideals. In Section \ref{sec.p-seqbei}, we shown that binomial edge ideals of trees and class of unicyclic graphs form a $p$-sequence. We study the Rees algebra of binomial edge ideals of trees in Sections \ref{sec.con}. In Section \ref{sec.counterex}, we present a counterexample to Conjecture \ref{conjlineartype} for a class of unicyclic graphs.  

\section{Preliminaries}

In this section, we recall basic definitions, notations, and results from the literature that are utilized in the subsequent sections.

\subsection{Basic notions from graph theory.} Let $G=(V(G)=[n],E(G))$ be a finite simple graph. For a subset $U \subseteq V(G)$, $G[U]$ denotes \textit{induced subgraph on $U$} of $G$, defined as $G[U]=(U,\{\{i,j\}\in E(G) \mid i,j \in U\})$. A graph $G$ is said to be \textit{path} if $E(G)=\{\{i,i+1\} \mid 1 \leq i \leq n-1\}\}$, and it is denoted by $P_n$. A graph $G$ is said to be  \textit{cycle} if $E(G)=\{ \{1,n\} \cup \{i,i+1\} \mid 1 \leq i \leq n-1\}\}$, for $n \geq 3$, and it is denoted by $C_n$. A graph $G$ is said to be \textit{forest} if it does not have a cycle and a \textit{tree} if it is connected.
It is said to be \textit{unicyclic} graph if it contains precisely one cycle as a subgraph.

\begin{notation} \label{not.g}
Let $G=(V(G),E(G))$ be a graph. For an edge $e$ and a vertex $v$ in $G$:
\begin{itemize}
    \item $G\setminus e$ is the graph on the vertex set $V(G)$ and edge set $E(G)\setminus {e}$.
    \item $e$ is called a \emph{bridge} if $c(G) < c(G \setminus e)$, where $c(G)$ is the number of components of $G$.
    \item $N_{G}(v) = \{u \in V(G) \mid \{u,v\} \in E(G)\}$ denotes the \emph{neighborhood} of $v$ in $G$.
    \item The \textit{degree} of a vertex $v \in V(G)$, denoted by $\deg_{G}(v)$, is defined as the number of edges incident to $v$. If $\deg_{G}(v)=1$, then $v$ is called a pendant vertex.
    \item An edge $e$ is called a pendant edge if one of the vertices is a pendant vertex.
    \item \cite[Definition 3.1]{MS} Let $e' =\{i,j\} \notin E(G)$ be an edge in $G \cup \{e'\}$. Then $G_{(e')}$ is the graph on vertex set $V(G)$ and edge set $$E(G_{(e')}) = E(G) \cup \left\{ \{k,l\} \mid k,l \in N_{G}(i) \textnormal{ or } k,l \in N_{G}(j) \right\}.$$
    %\item $V_{E(G)}=\{u\in V(G) \mid u \in e \subset E(G)\}$.
\end{itemize}
   
\end{notation}

A rooted tree $(G,x_0)$ is a directed tree in which edges are implicitly directed away from the root pendant vertex $x_0$. A level of a vertex $u \in V(G)$ in a rooted tree $(G,x_0)$ is the length of the path $P$ from $x_0$ to $u$. A level of an edge $e=\{x_i,x_j\}$ in a rooted tree $(G,x_0)$ is the length of the path $P$ from $x_0$ to $x_j$ such that $x_i \in V(P)$.

\subsection{Basic results from commutative algebra.}

Let $G$ denote a finite simple graph on vertex set $[n]$. Define $S=\mathbb{K}[x_1,\ldots,x_n,y_1,\ldots,y_n]$ as a polynomial ring over the field $\mathbb{K}$, and let $J_G$ be an ideal generated by edge binomials that are derived from $G$. Throughout this article, the symbols $S$ and $J_G$ are exclusively used to refer to the polynomial ring and the binomial edge ideal, respectively. The term ``sequence of elements in a ring" refers to a list of elements within that ring.

\vspace{2mm}

The subsequent results describe the colon ideal operation of an edge binomial on the binomial edge ideal, which we use repeatedly throughout this article.

\begin{remark}\cite[Lemma 4.1]{JAR2021} \cite[Corollary 2.2]{HH} \label{rem.colonp} Let $G$ be a simple graph and $f$ be an element of $S$. Then $J_G:f=J_G:f^n$ for all $n \geq 2$.
\end{remark} 

\begin{remark}\label{Rem.MCI} 
\cite[Theorem 3.4]{MS} Let $G$ be a simple graph. Let $e=\{i,j\} \notin E(G)$ be a bridge in $G \cup \{e\}$. Then $J_G: f_e = J_{G_e}$. 
\end{remark}

\begin{remark}\label{Rem.MCIunicy} \cite[Theorem 3.7]{MS}
Let $G$ be a simple graph and $e =\{i,j\} \notin E(G)$ be an edge. Then 
\[
J_G : f_e = J_{G_e} + (g_{P,t} \mid \quad P:i,i_1,\ldots,i_s,j \text{ is a path between } i,j \text{ and } 0 \leq t \leq s),
\] 
where $g_{P,0}=x_{i_1}\ldots x_{i_s}$ and for each $1\leq t \leq s, g_{P,t}=y_{i_1}\ldots y_{i_t}x_{i_{t+1}}\ldots x_{i_s}$.
\end{remark}

\begin{remark} \cite[Corollary 2.5]{SK} \label{rem.syz}
    For any bipartite graph $G$, one has $\beta_{i,i+1}(S/J_G) = 0$, for all $i$.
\end{remark}

\begin{definition}\cite[Definition 1.1]{H1982} \label{Def-d-sequence} Let $A$ be a commutative ring. Set $z_0 = 0$. A sequence of elements $z_1,\ldots,z_m$ in $A$ is said to be a $d$-sequence if it satisfies the conditions: 
\begin{enumerate}
    \item $z_1,\ldots,z_n$ is a minimal system of generators of the ideal $(z_1,\ldots,z_n)$;
    \item $( ( z_0,\ldots,z_i ) : z_{i+1}z_j) = (( z_0,\ldots,z_i ) : z_j)$ for all $0 \leq i \leq n-1$, and $j \geq i+1$.
\end{enumerate} 
     If $z_1,\ldots,z_m$ is a $d$-sequences in any order then the sequence is said to be  unconditional $d$-sequence.
\end{definition}

\begin{definition} \cite[Definition 1.1]{R1994} \label{def.relation}
    A \textit{relation} on an ordered sequence of elements $z_1,\ldots, z_n$ of a commutative ring $A$ is a form of positive degree in the graded polynomial ring $R=A[X_i \mid 1 \leq i \leq n]$, where elements of $A$ have degree $0$ and each $X_i$ has degree $1$, that yields $0$ when evaluated at $z_1,\ldots, z_n$. Let $\Lambda_j$ denote the ideal of $R$ generated by the relations on $z_1,\ldots, z_n$ of degree at most $j$. Let $r$ be the least positive integer such that $\Lambda_r =\cup_{j\geq r} \Lambda_j$. The invariant $r$ is called the \textit{relation type} of the ideal $I=(z_1,\ldots, z_n)$. An ideal of relation type $1$ is said to be of \textit{linear type}.
\end{definition}

Throughout this paper, we reserve $R$ for the polynomial ring over $S=\mathbb{K}[x_1,\ldots,x_n,y_1,\ldots,y_n]$ unless otherwise specified.

\begin{definition} \cite[p. 257]{C1985}
    Let $A$ be a commutative ring. An ordered sequence of elements $z_1,\ldots, z_n$ of a ring $A$ is said to be a \textit{sequence of linear type} if each of the ideal $(z_1,\ldots, z_i)$ is of linear type for all $1\leq i \leq n$.
\end{definition}

\section{\textit{p}-sequence} \label{sec.p-seq}

In this section, we introduce a notion called $p$-sequence. 

\begin{definition} \label{p-sequence}
    Let $A$ be a commutative ring. Set $z_0=0$.  A sequence of elements $z_1,\ldots,z_n$ in $A$ is said to be $p$-sequence if it satisfies the conditions: 
\begin{enumerate}
    \item \label{p-sequ.1} $z_1,\ldots,z_n$ is a minimal system of generators of the ideal $(z_1,\ldots,z_n)$;
    \item \label{p-sequ.2} $( z_0,\ldots,z_i ) : z_{i_1}z_{i_2} = ( z_0,\ldots,z_i ) : z_{i_1}$ for all $0 \leq i \leq n-1$, and $i < i_1 \leq i_2 \leq n$.

\end{enumerate} 

An ideal is said to be a $p$-sequence if a minimal generating set of the ideal forms a $p$-sequence.
\end{definition}

Clearly, unconditional $d$-sequences are $ p$-sequences. Examples of a class of unconditional $d$-sequences were given in \cite{H1982}. We list some examples of $p$-sequences here.

\begin{example}
    Any regular sequence is a $p$-sequence. 
\end{example}

\begin{proposition} \label{prop.p-subseq}
    Let $A$ be a commutative ring. Suppose the sequence of elements $z_1,\ldots,z_n$ in $A$ forms a $p$-sequence. Then $z_{i+1},\ldots,z_n$ in $A/I_i$ forms a $p$-sequence, where $I_i=(z_1,\ldots,z_{i})$, for all $i$.
\end{proposition}
\begin{proof}
    It follows from Definition \ref{p-sequence}.
\end{proof}

In the subsequent examples, we show that up to the sequence level, $d$-sequences and $p$-sequences are distinct notions. 
\begin{example}
    Let $G$ be a graph with vertex set $V(G)=[8]$ and edge set $$E(G)=\{\{1,2\},\{1,3\},\{1,4\},\{1,5\},\{5,6\},\{5,7\},\{5,8\}\}.$$ Then from Theorem \ref{thm.pseqtree}, it follows that there exists a sequence of edge binomials of $G$ that form a $p$-sequence.  However, the edge binomials of $G$ does not form a $d$-sequence under any permutation . This is due to the fact that for any graph $G$ having at least two vertices with degrees greater than or equal to $4$, edge binomials of $G$ can not form a $d$-sequence (cf. \cite[Corollary 3.2]{AN}).
\end{example} 
\begin{example}  \label{rem.pseqdseq}
    Consider a sequence of  elements $\{x_1x_3x_4x_5,x_1^2x_2x_6,x_1^2x_2^2x_3x_5\}$ in $\mathbb{K}[x_1,\ldots,x_6]$. In the specified order, this sequence readily forms a $d$-sequence. However, no permutation of this sequence constitutes a $p$-sequence. This assertion can be verified using Macaulay2 \cite{M2}.
\end{example}

\begin{definition}
    Let $A=\mathbb{K}[x_1,\ldots,x_n]$ be a polynomial ring over a field $\mathbb{K}$. Any product of variables of the form $x^\alpha = x_1^{a_1}\cdots x_n^{a_n}$ is called monomial, where $\alpha \in \mathbb{N}^n$.
\end{definition}

\begin{proposition}\label{prop.p-mon}
Let \(z_1,\ldots,z_n\) be an ordered sequence of monomials in \(\mathbb{K}[x_1,\ldots,x_n]\). Suppose that \(z_i \nmid z_j\) for all \(i \neq j\), and that the following conditions hold:
\begin{equation}\label{cond1}
\gcd(z_i,z_{i_2}) \mid z_{i_1},
\qquad
1 \leq i < i_1 \leq i_2 \leq n,
\end{equation}
\begin{equation}\label{cond2}
\gcd(z_i,z_j^2) = \gcd(z_i,z_j),
\qquad
1 \leq i < j \leq n.
\end{equation}
Then \(z_1,\ldots,z_n\) is a \(p\)-sequence.
\end{proposition}

\begin{proof}
Since $z_i\nmid z_j$ for all $i\neq j$, the sequence $z_1,\ldots,z_n$ is a minimal system of monomial generators. To prove that \(z_1,\ldots,z_n\) is a \(p\)-sequence, it suffices to show that
\[
\gcd(z_i, z_{i_1}z_{i_2}) = \gcd(z_i, z_{i_1})
\quad \text{for all } 1 \leq i < i_1 \leq i_2 \leq n.
\]

Fix \(i < i_1 \leq i_2\). By \eqref{cond1}, \(\gcd(z_i, z_{i_2}) \mid z_{i_1}\); hence, it follows that \(\gcd(z_i, z_{i_1}z_{i_2}) \mid \gcd(z_i, z_{i_1}^2)\). Using \eqref{cond2}, we obtain $\gcd(z_i, z_{i_1}^2) = \gcd(z_i, z_{i_1})$,
and hence
$\gcd(z_i, z_{i_1}z_{i_2}) \mid \gcd(z_i, z_{i_1})$.
On the other hand, it is clear that
$\gcd(z_i, z_{i_1}) \mid \gcd(z_i, z_{i_1}z_{i_2})$.
Therefore,
$\gcd(z_i, z_{i_1}z_{i_2}) = \gcd(z_i, z_{i_1})$,
for all \(1 \leq i < i_1 \leq i_2 \leq n\). This proves that \(z_1,\ldots,z_n\) is a \(p\)-sequence.
\end{proof}

\begin{remark}
    In the context of edge ideals, Proposition~\ref{prop.p-mon} yields an ordering of the generators under which the edge ideals of matching and caterpillar graphs are generated by a \(p\)-sequence.  It is known from Villarreal~\cite{V1990} that the edge ideals of these graphs are of linear type.
    %These graphs evidently belong to the subset of ideals characterized as ideals of linear type by Villareal \cite{V1990}. Thus, it prompts the natural question.
\end{remark}

This naturally leads to the following question.

\begin{question}\label{qus.linear}
Let \(A\) be a commutative ring, and let \(z_1,\ldots,z_n\) in  \(A\) be a \(p\)-sequence. Is the ideal \(I=(z_1,\ldots,z_n)\) of linear type?
\end{question}

The following example demonstrates that a $p$-sequence fails to satisfy one of the properties of a $d$-sequence.
\begin{example}
    Suppose the sequence of elements $z_1,\ldots,z_n$ in $S$ form a $d$-sequence. Then, it satisfies  $$z_i \cap ( z_1,\ldots,z_{i-1})^s(z_1,\ldots,\hat{z_i},\ldots,z_n)\subseteq z_i(z_1,\ldots,z_{i-1})^{s-1}(z_1,\ldots,\hat{z_i},\ldots,z_n)$$ (see \cite[Proposition 3.1]{H1980}). However, this property does not hold for $p$-sequences. For instance, consider a graph $G$ with edge set $$E(G)=\{\{1,2\},\{2,3\},\{2,4\},\{4,5\},\{4,6\},\{3,5\}\}.$$ From Theorem \ref{thm.pseqcycle}, it follows that the ideal $J_G$ forms a $p$-sequence. Then, one can verify that $$f_{24}\cap (f_{12},f_{23})^2(f_{12},f_{23},f_{45},f_{46},f_{35}) \not\subset f_{24}(f_{12},f_{23})(f_{12},f_{23},f_{45},f_{46},f_{35}).$$
\end{example}

\section{Edge binomials of trees and unicyclic graphs form a $p$-sequence} \label{sec.p-seqbei}

In this section, we demonstrate that the edge binomials of trees and a class of unicyclic graphs form a $p$-sequence.

\subsection{Ordering on edges of a tree} \label{oderseqtree} Let $G$ be a tree. Any tree can be represented as a rooted tree structure with a pendant vertex as a root vertex, denoted by $(G,x_0)$. For each level \(m\), let 
\[
e_{m,1}, e_{m,2}, \ldots, e_{m,k_m}
\]
be the edges at level \(m\), listed from left to right according to a fixed planar embedding of \(G\). Since every tree admits a planar embedding, such an ordering is well defined. We ordered the edge binomials of a tree in the following order. Let $(G,x_0)$ be a rooted tree with maximum level $m$. Then we arrange the edges  by $$e_{1,1},e_{2,1},\ldots,e_{2,k_2},\ldots,e_{m-1,1},\ldots, e_{m-1,k_{m-1}}, e_{m,1},\ldots,e_{m,k_m}.$$

\begin{example}
    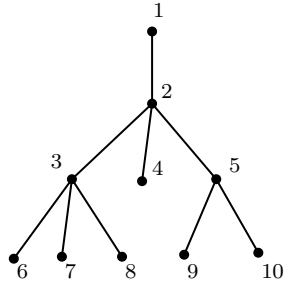
\begin{figure}[h]
    \centering

\tikzset{every picture/.style={line width=0.75pt}} %set default line width to 0.75pt        

\begin{tikzpicture}[x=0.75pt,y=0.75pt,yscale=-1,xscale=1]
%uncomment if require: \path (0,2583); %set diagram left start at 0, and has height of 2583

%Straight Lines [id:da7216796857582822] 
\draw    (322,2134.58) -- (322,2170.88) ;
%Straight Lines [id:da6623309344823036] 
\draw    (322,2170.88) -- (317,2209.88) ;
%Straight Lines [id:da4392242490575906] 
\draw    (322,2170.88) -- (282,2208.88) ;
%Straight Lines [id:da8456341042665158] 
\draw    (322,2170.88) -- (354,2208.88) ;
%Straight Lines [id:da0696179057810431] 
\draw    (282,2208.88) -- (277,2247.88) ;
%Straight Lines [id:da2565237837748221] 
\draw    (282,2208.88) -- (307,2247.88) ;
%Straight Lines [id:da3294824072435606] 
\draw    (282,2208.88) -- (253,2248.88) ;
%Straight Lines [id:da4113883095362221] 
\draw    (354,2208.88) -- (338,2246.88) ;
%Straight Lines [id:da05832457593090201] 
\draw    (354,2208.88) -- (375,2245.88) ;

% Text Node
\draw (321,2118.78) node [anchor=north west][inner sep=0.75pt]  [font=\scriptsize]  {$1$};
% Text Node
\draw (325,2159.78) node [anchor=north west][inner sep=0.75pt]  [font=\scriptsize]  {$2$};
% Text Node
\draw (270,2194.78) node [anchor=north west][inner sep=0.75pt]  [font=\scriptsize]  {$3$};
% Text Node
\draw (320.5,2198.38) node [anchor=north west][inner sep=0.75pt]  [font=\scriptsize]  {$4$};
% Text Node
\draw (359,2196.78) node [anchor=north west][inner sep=0.75pt]  [font=\scriptsize]  {$5$};
% Text Node
\draw (253,2250) node [anchor=north west][inner sep=0.75pt]  [font=\scriptsize]  {$6$};
% Text Node
\draw (277,2250) node [anchor=north west][inner sep=0.75pt]  [font=\scriptsize]  {$7$};
% Text Node
\draw (307,2250) node [anchor=north west][inner sep=0.75pt]  [font=\scriptsize]  {$8$};
% Text Node
\draw (338,2250) node [anchor=north west][inner sep=0.75pt]  [font=\scriptsize]  {$9$};
% Text Node
\draw (375,2250) node [anchor=north west][inner sep=0.75pt]  [font=\scriptsize]  {$10$};

%\draw (300,2145) node [anchor=north west][inner sep=0.75pt]  [font=\scriptsize]  {$e_{1,1}$};

%\draw (285,2180) node [anchor=north west][inner sep=0.75pt]  [font=\scriptsize]  {$e_{2,1}$};

%\draw (300,2190) node [anchor=north west][inner sep=0.75pt]  [font=\scriptsize]  {$e_{2,2}$};

\filldraw[black] (322,2134.58) circle (1.5pt) ;
\filldraw[black] (322,2170.88) circle (1.5pt) ;
\filldraw[black] (317,2209.88) circle (1.5pt) ;
\filldraw[black] (282,2208.88) circle (1.5pt) ;
\filldraw[black] (354,2208.88)  circle (1.5pt) ;
\filldraw[black] (277,2247.88) circle (1.5pt) ;
\filldraw[black] (307,2247.88) circle (1.5pt) ;
\filldraw[black] (253,2248.88) circle (1.5pt) ;
\filldraw[black] (338,2246.88)  circle (1.5pt) ;
\filldraw[black] (375,2245.88)  circle (1.5pt) ;
\end{tikzpicture}

\caption{The rooted tree $(G,1)$}
    \label{fig:rootedtree}
    \end{figure}

    Let $(G,1)$ be a rooted tree as given in Figure \ref{fig:rootedtree}. Then, we arrange the edges of $(G,1)$ in the following order $e_{1,1}=\{1,2\},e_{2,1}=\{2,3\},e_{2,2}=\{2,4\},e_{2,3}=\{2,5\},e_{3,1}=\{3,6\},e_{3,2}=\{3,7\},e_{3,3}=\{3,8\},e_{3,4}=\{5,9\},e_{3,5}=\{5,10\}$.

\end{example}

\subsection{Ordering on edges of a unicyclic graph} \label{odersequnicycle} 
Let $(T,1)$ be a rooted tree with the maximal level of a vertex is $m$. Let $G$ be a unicyclic graph obtained by adding an edge $e=\{u,v\}$ between two pendant vertices of tree $T$. Represent the tree $T$ by rooted tree structure as in \ref{oderseqtree} such that one of the vertices of $e$ is left most in level $m$ vertex. Then we arrange the edges of $G=(T,1) \cup e$ by  $e_{1,1},e_{2,1},\ldots,e_{2,k_2},\ldots,e_{m-1,1},\ldots, e_{m-1,k_{m-1}}, e_{m,1},\ldots,e_{m,k_m},e$.

\begin{example}
    Let $G_1 = (G,1)\cup \{8,10\}$ (resp. $G_1 = (G,1)\cup \{4,10\}$) be a unicyclic graph, where $(G,1)$ is a rooted tree as given in Figure \ref{fig:rootedtree}. Then, we arrange the edges of $G_1$ in the following order $e_{1,1}=\{1,2\},e_{2,1}=\{2,3\},e_{2,2}=\{2,4\},e_{2,3}=\{2,5\},e_{3,1}=\{3,6\},e_{3,2}=\{3,7\},e_{3,3}=\{3,8\},e_{3,4}=\{5,9\},e_{3,5}=\{5,10\},e=\{8(\text{resp. } 4),10\}$.
\end{example}

\begin{theorem}\label{thm.pseqtree}
Let \(G\) be a tree on \([n+1]\). Then the binomial edge ideal of \(G\) forms a \(p\)-sequence ideal.
\end{theorem}

\begin{proof}
We show that the ordered sequence of edge binomials associated to \(G\), arranged as in Subsection~\ref{oderseqtree}, forms a \(p\)-sequence. Let \(d_1,\ldots,d_n\) be the edge binomials of \(G\) in the prescribed order, and set \(d_0=0\). For each \(i\), let \(e_i\) denote the edge corresponding to \(d_i\). It suffices to verify that
\[
(d_0,\ldots,d_i):d_{i_1}d_{i_2}
=
\big((d_0,\ldots,d_i):d_{i_1}\big):d_{i_2}
=
(d_0,\ldots,d_i):d_{i_1},
\]
for all \(0 \le i \le n-1\) and \(i < i_1 \le i_2 \le n\).

Fix \(i < n\), and let \(H\) be the subgraph of \(G\) with edge set \(\{e_1,\ldots,e_i\}\). By the ordering of the edges (cf. Subsection~\ref{oderseqtree}), edge \(e_{i_1}\) with \(i_1>i\) intersects \(H\) in at most one vertex. By Remark~\ref{Rem.MCI}, we have
\[
(d_0,\ldots,d_i):d_{i_1} = J_{H_{(e_{i_1})}}.
\]

Let \(e_{i_2}=\{u,v\}\). We analyze the colon ideal \(J_{H_{(e_{i_1})}}:d_{i_2}\) according to how \(e_{i_2}\) intersects \(H\) and \(e_{i_1}\).

\medskip

\noindent
\textit{Case 1:} \(e_{i_2} \cap H = \emptyset\).  
In this case, Remark~\ref{Rem.MCI} yields
\[
J_{H_{(e_{i_1})}}:d_{i_2} = J_{H_{(e_{i_1})}}.
\]

\medskip

\noindent
\textit{Case 2:} \(e_{i_2} \cap H \neq \emptyset\) and \(e_{i_2} \cap e_{i_1} = \emptyset\).  
By the edge ordering, exactly one vertex of \(e_{i_2}\) lies in \(H\). Without loss of generality, assume \(u \in V(H)\) and \(v \notin V(H)\). Then in the graph \(H_{(e_{i_1})} \cup e_{i_2}\), we have
\[
|N(u)| = 2 \quad \text{and} \quad |N(v)| = 1.
\]
Applying Remark~\ref{Rem.MCI}, it follows that
\[
J_{H_{(e_{i_1})}}:d_{i_2} = J_{H_{(e_{i_1})}}.
\]

\medskip

\noindent
\textit{Case 3:} \(e_{i_2} \cap H \neq \emptyset\) and \(e_{i_2} \cap e_{i_1} \neq \emptyset\).  
In this case, either \(e_{i_2}=e_{i_1}\) or the two edges share a common vertex. In both situations, the neighborhoods of the vertices of \(e_{i_2}\) in \(H\) coincide with those of \(e_{i_1}\). Hence, by Remark~\ref{Rem.MCI},
\[
J_{H_{(e_{i_1})}}:d_{i_2} = J_{H_{(e_{i_1})}}.
\]

\medskip

In all cases, we obtain
\[
\big((d_0,\ldots,d_i):d_{i_1}\big):d_{i_2}
=
(d_0,\ldots,d_i):d_{i_1}.
\]
Therefore, \(d_{i_2}\) is a non-zero divisor modulo \((d_0,\ldots,d_i):d_{i_1}\), and the sequence \(d_1,\ldots,d_n\) is a \(p\)-sequence. This completes the proof.
\end{proof}

\begin{theorem}\label{thm.pseqcycle}
Let \(T\) be a tree on \([n]\), and let \(G\) be the graph obtained by adding an edge between two pendant vertices of \(T\). Then the edge binomials of \(G\) form a \(p\)-sequence.
\end{theorem}

\begin{proof}
Let \(d_1,\ldots,d_n\) be the edge binomials of \(G\) ordered as in Subsection~\ref{odersequnicycle}, and set \(d_0=0\). For each \(i\), let \(e_i\) denote the edge corresponding to \(d_i\). We show that the sequence \(d_1,\ldots,d_n\) is a \(p\)-sequence.

By Theorem~\ref{thm.pseqtree}, we already know that
\[
(d_0,\ldots,d_i):d_{i_1}d_{i_2}
=
(d_0,\ldots,d_i):d_{i_1}
\]
for all \(0 \le i \le n-2\) and \(i < i_1 \le i_2 \le n-1\). Thus, it remains to verify that
\[
(d_0,\ldots,d_i):d_{i_1}d_n
=
(d_0,\ldots,d_i):d_{i_1}
\]
for all \(0 \le i \le n-1\) and \(i < i_1 \le n\).

Fix \(i\) and \(i<i_1\le n\), and let \(H\) be the subgraph of \(G\) with edge set \(\{e_1,\ldots,e_i\}\). By construction of the ordering, for \(i<n\), both edges \(e_{i_1}=\{r,s\}\) and \(e_n=\{u,v\}\) are leaves when attached to the graph \(H\). Without loss of generality, assume that \(s\) and \(v\) are pendant vertices in \(H \cup e_{i_1}\) and \(H \cup e_n\), respectively.

\medskip

\noindent
\textit{Case 1:} \(i_1 < n\).  
By Remark~\ref{Rem.MCI}, we have
\[
(d_0,\ldots,d_i):d_{i_1} = J_{H_{(e_{i_1})}}.
\]
We analyze the colon ideal \(J_{H_{(e_{i_1})}}:d_n\).

If \(e_n \cap H = \emptyset\), then we have $J_{H_{(e_{i_1})}}:d_n = J_{H_{(e_{i_1})}}$
by Remark~\ref{Rem.MCI}.

If \(e_n \cap H \neq \emptyset\), then there exists \(e_j\), with \(j \le i\), such that \(e_n \cap e_j \neq \emptyset\). Suppose \(u \in V(H)\). If \(\deg_{H_{(e_{i_1})} \cup e_n}(u)=2\), then again Remark~\ref{Rem.MCI} yields
\[
J_{H_{(e_{i_1})}}:d_n = J_{H_{(e_{i_1})}}.
\]
If \(\deg_{H_{(e_{i_1})} \cup e_n}(u)>2\), then by the ordering of edges, the neighborhoods satisfy
$N_{H}[r]=N_{H_{(e_{i_1})}}[u]$,
and hence, by Remark~\ref{Rem.MCI}, we again obtain
\[
J_{H_{(e_{i_1})}}:d_n = J_{H_{(e_{i_1})}}.
\]

Thus, in all cases,
\[
(d_0,\ldots,d_i):d_{i_1}d_n
=
(d_0,\ldots,d_i):d_{i_1}.
\]

\medskip

\noindent
\textit{Case 2:} \(i_1 = n\).  
In this case, the equality
\[
(d_0,\ldots,d_i):d_n^2
=
(d_0,\ldots,d_i):d_n
\]
follows directly from Remark~\ref{rem.colonp}.

\medskip

Combining the above, we conclude that the sequence \(d_1,\ldots,d_n\) satisfies the defining property of a \(p\)-sequence. Hence, the edge binomials of \(G\) form a \(p\)-sequence.
\end{proof}

\section{Proof of Conjecture \ref{conjlineartype} for trees} \label{sec.con}
In this section, we aim to prove Conjecture \ref{conjlineartype} for all trees. We begin by establishing the necessary notations.

\begin{notation}
    Let $G=(V(G),E(G))$ be a graph with $\lvert E(G) \rvert = m$. We arrange the edge of $G$ in an ordered sequence $e_1,\ldots,e_m$. We denote $G_i$ as the graph with vertex set $V(G_i)=V(G)$ and edge set $E(G_i)=\{e_1,\ldots,e_i\}$, for $1\leq i \leq m$. The graph $G_{i_{(e_{i+1})}}$ is defined as in Notation \ref{not.g}.
\end{notation}

\begin{lemma}\label{lem.colon.bi.ideal}
Let \(G\) be a tree on vertex set \([n+1]\), and let \(e_1,\ldots,e_n\) be an ordered sequence of edges of \(G\) as described in Subsection~\ref{oderseqtree}. For each \(i=1,\ldots,n-1\), let \(T_i\) be the subgraph of \(G\) with vertex set \(V(T_i)=V(G)\) and edge set \(E(T_i)=E(G)\setminus E(G_i)\). Then
\begin{equation}\label{eq.colon}
J_{G_{i_{(e_{i+1})}}}:J_{T_i} = J_{G_{i_{(e_{i+1})}}}
\end{equation}
for all \(i=1,\ldots,n-1\). Moreover, for any nonzero \(f \in J_{T_i}\), one has
\[
J_{G_{i_{(e_{i+1})}}}:f = J_{G_{i_{(e_{i+1})}}}.
\]
\end{lemma}

\begin{proof}
Fix \(i\in\{1,\ldots,n-1\}\). Then \(E(T_i)=\{e_{i+1},\ldots,e_n\}\), and hence
\[
J_{G_{i_{(e_{i+1})}}}:J_{T_i}
=
\bigcap_{j=i+1}^{n} \big( J_{G_{i_{(e_{i+1})}}} : f_{e_j} \big).
\]

We show that \(J_{G_{i_{(e_{i+1})}}} : f_{e_j} = J_{G_{i_{(e_{i+1})}}}\) for each \(j \ge i+1\).  
Let \(e=\{u,v\}\) be an edge of \(T_i\). Since \(G\) is a tree and by the chosen ordering, every such edge \(e\) is a bridge in the graph \(G_i \cup e\).

We consider the possible intersections of \(e\) with \(G_i\) and \(e_{i+1}\).

\medskip

\noindent
\textit{Case 1:} \(e \cap G_i = \emptyset\).  
Then \(e\) is disjoint from \(G_i\), and by Remark~\ref{Rem.MCI},
\[
J_{G_{i_{(e_{i+1})}}} : f_e = J_{G_{i_{(e_{i+1})}}}.
\]

\medskip

\noindent
\textit{Case 2:} \(e \cap G_i \neq \emptyset\) and \(e \cap e_{i+1} = \emptyset\).  
In this case, exactly one vertex of \(e\) lies in \(G_i\). Without loss of generality, assume \(u \in V(G_i)\). Then
\[
\deg_{G_i \cup e}(u)=2 \quad \text{and} \quad \deg_{G_i \cup e}(v)=1.
\]
Hence, by Remark~\ref{Rem.MCI},
\[
J_{G_{i_{(e_{i+1})}}} : f_e = J_{G_{i_{(e_{i+1})}}}.
\]

\medskip

\noindent
\textit{Case 3:} \(e \cap G_i \neq \emptyset\) and \(e \cap e_{i+1} \neq \emptyset\).  
Then either \(e=e_{i+1}\) or \(e\) and \(e_{i+1}\) share a common vertex. In both cases, the neighborhoods of the common vertices of \(e\) and \(e_{i+1}\) in \(G_i\) coincide. Thus, by Remark~\ref{Rem.MCI},
\[
J_{G_{i_{(e_{i+1})}}} : f_e = J_{G_{i_{(e_{i+1})}}}.
\]

\medskip

Therefore, $J_{G_{i_{(e_{i+1})}}}:J_{T_i} = \bigcap_{j=i+1}^{n} J_{G_{i_{(e_{i+1})}}} = J_{G_{i_{(e_{i+1})}}}$, which proves \eqref{eq.colon}.

\medskip

For the second assertion, the inclusion
\[
J_{G_{i_{(e_{i+1})}}} \subseteq J_{G_{i_{(e_{i+1})}}}:f
\]
is immediate. For the reverse inclusion, let \(a \notin J_{G_{i_{(e_{i+1})}}}\) and suppose \(af \in J_{G_{i_{(e_{i+1})}}}\) for some nonzero \(f \in J_{T_i}\). Then, by Remark~\ref{Rem.MCIunicy}, there exists a path in \(T_i\) connecting two vertices \(\alpha,\beta \in V(G_{i_{(e_{i+1})}})\). This contradicts the fact that \(G\) is a tree. Hence
\[
J_{G_{i_{(e_{i+1})}}}:f \subseteq J_{G_{i_{(e_{i+1})}}},
\]
and equality follows.
\end{proof}

\begin{notation}
Let \(z_1,\ldots,z_t\) be an ordered sequence of elements in \(S\). Consider the homomorphism
\[
\phi : R \longrightarrow S, \qquad X_i \mapsto z_i \ \text{for } i=1,\ldots,t.
\]
For each \(i\), denote by \(I_i\) the ideal of \(R\) generated by \(X_1,\ldots,X_i\), and by \(\overline{I}_i\) the ideal of \(S\) generated by \(z_1,\ldots,z_i\).

For a polynomial \(G(X_1,\ldots,X_t) \in R\), we write
\[
\overline{G} := G(z_1,\ldots,z_t) \in S
\]
for its image under \(\phi\).
\end{notation}

\begin{theorem} \label{thm.conj}
    Let $G$ be a tree on $[n+1]$. Then the ideal $J_G$ is of linear type. 
\end{theorem}

\begin{proof}
     Let $f_{e_1},\ldots,f_{e_n}$ be the ordered sequence of edge binomials of $G$ as described in Subsection \ref{oderseqtree}. By Theorem \ref{thm.pseqtree} this sequence forms a $p$-sequence. We claim that for any polynomial $F(X_1,\ldots,X_n)$ in $R$, if evaluated at $(f_{e_1},\ldots,f_{e_n})$ results in $0$, then $F$ is of linear relation type. We prove this by strong induction on $\deg(F)$. 
     
     If $\deg(F)=1$, there is nothing to prove. Assume that the statement holds for all $\deg(F) < d$. For $\deg(F)=d$, let $F=\sum_{\lvert \alpha \rvert \leq d}a_{\alpha}X^{\alpha}$ be a polynomial such that evaluated at $f_{e_1},\ldots,f_{e_n}$ yields $0$, where $\alpha \in \mathbb{N}^n$. 
     
     \vspace{2mm}

     Set $F_0=F$ and $I_0=0$. The proof proceeding in following steps as $k$ runs from $1$ to $n$. Here, we define $F_k$ recursively in Step $k$. 
     
     \vspace{2mm}

    \textbf{Step $k$.}
    
    Considering, $F_{k-1}$ mod $I_{n-k}$  have the following cases:

    \vspace{2mm}
    
    Case I. If $F_{k-1}(X_1,\ldots,X_n) = aX_{n-(k-1)}^{\ell_1}\cdots X_n^{\ell_k} \mod I_{n-k}$, where $\lvert \ell_1 \rvert + \cdots + \lvert \ell_k \rvert \leq d$. Since $\overline{F}=0$ it follows that $af_{e_{n-(k-1)}}^{\ell_1}\cdots f_{e_{n}}^{\ell_k} \in \overline{I}_{n-k}$. From Proposition \ref{prop.p-subseq} and $p$-sequence property, it follows that  
    $$a \in \overline{I}_{n-k}:f_{e_{n-(k-1)}}^{\ell_1}\cdots f_{e_{n}}^{\ell_k} = \overline{I}_{n-k}:f_{e_{n-(k-1)}}.$$ 
    Hence $af_{e_{n-(k-1)}}-\sum_{i=1}^{n-k}b_if_{e_i} = 0$. Set $H_k=aX_{n-(k-1)}-\sum_{i=1}^{n-k}b_iX_i$ and $$F_k=F_{k-1}-H_k(X_{n-(k-1)}^{\ell_1-1}\cdots X_n^{\ell_k}).$$ Then $F_k \in (X_1,\ldots,X_{n-k})$ if $k<n$. If $k=n$ then $F_n$ has a linear relation.
    
    \vspace{2mm}
    
    Case II. Suppose $F_{k-1}(X_1,\ldots,X_n) = aX_{n-(k-1)}H \mod I_{n-k}$, where $H \in S[X_n,\ldots,X_{n-(k-1)}]$, $\deg(H)<d$ and $H$ is not a monomial. Since $\overline{F}=0$ we have $\overline{F}_{k-1} = af_{e_{n-(k-1)}}\overline{H}$ mod $\overline{I}_{n-k} = 0 $ mod $\overline{I}_{n-k}$. This implies that $a\overline{H} \in \overline{I}_{n-k}:f_{e_{n-(k-1)}}$.  Note that $\overline{I}_{n-k}=J_{G_{n-k}}$ and $e_{n-(k-1)}$ is a bridge in $G_{n-k} \cup e_{n-(k-1)}$. It follows from Remark \ref{Rem.MCI} that $a\overline{H} \in J_{G_{n-k_{(e_{n-(k-1)})}}}$. By Lemma \ref{lem.colon.bi.ideal},  $a \in J_{G_{n-k_{(e_{n-(k-1)})}}}: \overline{H}= J_{G_{n-k_{(e_{n-(k-1)})}}}$. Therefore, either $a\overline{H}=0$ or every monomial that survive in $a\overline{H}$ belong to $J_{G_{n-k_{(e_{n-(k-1)})}}}$.

    \vspace{2mm}
    
    Subcase A. If $a\overline{H}=0$, then by induction there is a linear relation $H'= \sum_{i=1}^nb_iX_i$ such that $\overline{H'}=0$, and $KH' = aH$ for some $K \in R$.
    Since $a\overline{H}=0$ and $\overline{F}=0$, we have $\overline{F_{k-1}}-a\overline{X_{n-(k-1)}H}=0$. Since $aH$ has linear relation, it suffices to show that  $F_{k-1}-a{X_{n-(k-1)}H}$ has linear relation. If $k<n$, set $$F_k=F_{k-1}-a{X_{n-(k-1)}H},$$ so that $F_k\in (X_1,\ldots,X_{n-k})$.  If $k=n$ then $F_n=a{X_{1}H}$ has a linear relation.
   
    \vspace{2mm}
    
    Subcase B. If $a\overline{H}\neq 0$ (i.e, there exists a monomial in the polynomial $aX_{n-(k-1)}H$ which survive on  $af_{e_{n-(k-1)}}\overline{H}$), then choose leading term of $aX_{n-(k-1)}H$ with respect to any fixed term order which survive in $af_{e_{n-(k-1)}}\overline{H}$, say $cX_{n-(k-1)}^{\ell_1}\cdots X_n^{\ell_k}$, where $\lvert \ell_1 \rvert + \cdots + \lvert \ell_k \rvert \leq d$ and $\ell_1 >0$. We show that each such surviving monomial can be linearized. The monomial $cf_{e_{n-(k-1)}}^{\ell_1}\cdots f_{e_{n}}^{\ell_k}$ belongs to $\overline{I}_{n-k}$, since $cf_{e_{n-(k-1)}} \in \overline{I}_{n-k}$ as $c \in J_{G_{n-k_{(e_{n-(k-1)})}}}$. Then from Proposition \ref{prop.p-subseq} and Definition \ref{p-sequence} it follows that $c \in \overline{I}_{n-k}:f_{e_{n-(k-1)}}^{\ell_1}\cdots f_{e_{n}}^{\ell_k} = \overline{I}_{n-k}:f_{e_{n-(k-1)}}$. Therefore one has $cf_{e_{n-(k-1)}}-\sum_{i=1}^{n-k}b_if_{e_i} = 0$. Set $H_{k,1}=cX_{n-(k-1)}-\sum_{i=1}^{n-k}b_iX_i$ and $$F_{k,1}= F_{k-1}-H_{k,1}(X_{n-(k-1)}^{\ell_1-1}\cdots X_n^{\ell_k}).$$ Moreover, $F_{k,1} \in (X_1,\ldots,X_{n-(k-1)})$.
    
    Repeat Case II for $F_{k,1}$. If it falls under either Case I or Subcase A, then proceed to step $k+1$. Otherwise, recursively define  $H_{k,j}$ and $F_{k,j}$, and continue applying Case II to $F_{k,j}$ for $j\geq 2$, until it eventually falls under either Case I or Subcase A. This recursive process terminates in finitely many steps, as $aX_{n-(k-1)}H$ is a polynomial. 
    
    \vspace{2mm}
    
    Case III. If $F_{k-1}(X_1,\ldots,X_n) = 0$ mod $I_{n-k}$ then set $F_k=F_{k-1}$ and clearly  $F_k\in (X_1,\ldots,X_{n-k})$

    \vspace{2mm}

    Since \(F_k \in (X_1,\ldots,X_{n-k})\) for all \(k=0,\ldots,n-1\), the subcase of the form $F_{k}(X_1,\ldots, X_n) = aX_{n-k}H_1 + H_2 \mod I_{n-(k-1)}$, where $H_1 \in S[X_{n-k},\ldots, X_n]$ and $H_2 \in S[X_{n-(k-1)},\ldots, X_n]$ do not occur.
    
    %Since $F_k\in (X_1,\ldots,X_{n-k})$ for all $k$, the subcase of the form $F_{k-1}(X_1,\ldots, X_n) = aX_{n-(k-1)}H_1 + H_2 \mod I_{n-k}$, where $H_1 \in S[X_{n-(k-1)},\ldots, X_n]$ and $H_2 \in S[X_{n-(k-2)},\ldots, X_n]$ is not possible.

    Observe that at each k$^{th}$ step, the claim  reduces to showing that $F_{k}$ is of linear type. Hence, by the end of Step $n$, one can conclude that $F$ is of linear type. 
\end{proof}

\begin{corollary}
     If $G$ be a tree then $J_G$ is a sequence of linear type.
\end{corollary}
\begin{proof}
    The proof follows from Theorem \ref{thm.pseqtree} and Theorem \ref{thm.conj}.
\end{proof}

In (\cite[Corollary 6.2]{B2018}), Bolognini et al. proved that if $G$ is a bipartite graph, then the binomial edge ideal $J_G$ coincides with the parity binomial edge ideal $\mathcal{I}_G$ (cf. \cite{TCT}) and Lovász–Saks–Schrijver ideal $L_G(2)$ (cf. \cite{HMSW}), use which we establish the linear type properties of these ideals for trees in the following remark.

\begin{remark}
    By using \cite[Corollary 6.2]{B2018}, one can obtain that if $G$ be a tree, then the associated parity binomial edge ideal and Lovász–Saks–Schrijver ideal are of linear type. 
\end{remark}

\begin{corollary}
Let \(z_1,\ldots,z_n\) be a \(p\)-sequence of elements in \(S\). Suppose that for each \(i=1,\ldots,n\), and for every polynomial \(H \in (z_1,\ldots,z_i)\), each term of \(H\) belongs to \((z_1,\ldots,z_i)\). Then the ideal \(I=(z_1,\ldots,z_n)\) is of linear type.
\end{corollary}
\begin{proof}
The statement follows from the argument in the proof of Theorem~\ref{thm.conj}. In particular, consider a polynomial of the form \(aX_{n-(k-1)}H\), where \(H \in S[X_n,\ldots,X_{n-(k-1)}]\). From the proof of Theorem~\ref{thm.conj}, one has $a\overline{H} \in \overline{I}_{n-k}$.

By assumption, every term of \(az_{n-(k-1)}\overline{H}\) belongs to \(\overline{I}_{n-k}\). It follows that each corresponding term in \(aX_{n-(k-1)}H\) can be reduced to linear relations. Consequently, the polynomial \(aX_{n-(k-1)}H\) has a linear relation.
\end{proof}

\subsection{Defining ideal of Rees algebra}
The subsequent result follows from Theorem \ref{thm.pseqtree} and \cite[Theorem 3.2]{JAR2021}, where we obtain a minimal generating set for the defining ideal of $\mathcal{R}(J_G)$ when $G$ is a tree.

\begin{corollary}
    Let $\sigma : S[X_{\{i,j\}} \mid \{i,j\} \in E(G)] \rightarrow \mathcal{R}(J_G)$ be the map defined by $\sigma(X_{\{i,j\}})= f_{ij}t$. If $G$ is a tree then the kernel of $\sigma$ (the defining ideal of $\mathcal{R}(J_G)$) is minimally generated by 
    $\{ f_{ij}X_{\{k,l\}}-f_{kl}X_{\{i,j\}} \mid \{i,j\} \neq \{k,l\} \in E(G)\} \cup \{(-1)^{p_A(j)}f_{kl}X_{\{i,j\}}+(-1)^{p_A(k)}f_{jl}X_{\{i,k\}}+(-1)^{p_A(l)}f_{jk}X_{\{i,l\}}$, where $A=\{i,j,k,l\}$ is claw $(K_{1,3})$ with center at $i$ and $p_A(i)= \lvert \{j \in A \mid j \leq i\} \rvert\}$.
\end{corollary}

\section{Counterexample} \label{sec.counterex}
In this section, we present a counterexample that refutes Conjecture \ref{conjlineartype} specifically for a certain family of unicyclic graphs. Initially, we define the notation necessary for our demonstration.

\begin{notation}
    We denote by $C_{n,k}$ the class of graphs obtained by attaching $k$ pendant edges to each vertex of the cycle graph $C_n$.
\end{notation}

\begin{example} \label{exa.c41}
    $C_{4,1}$ be a graph with vertex set $V(C_{4,1})=[8]$ and edge set $E(C_{4,1})=\{\{1,2\},\{2,3\}$, $\{3,4\},\{1,4\}$, $\{1,5\},\{2,6\},\{3,7\},\{4,8\}\}$.
\end{example}

\begin{proposition} \label{prop.counter}
    If $C_{4,1}$ is an induced subgraph of $G$, then $J_G$ is not a linear type.
\end{proposition}
\begin{proof}
    Let $\phi: S[X_{ij} \mid \{i,j\} \in E(G)] \rightarrow S$ be a map given by $\phi(X_{ij})=f_{ij}$. Let $\Lambda$ be a kernel of $\phi$ (see, Definition \ref{def.relation}). Suppose  $C_{4,1}$ is an induced subgraph of $G$. Without loss of generality, we may assume that $C_{4,1}$ has the same vertex set and edge set as given in Example \ref{exa.c41}. Then, one can verify that
    \begin{equation} \label{eq.qudratic}
        \begin{split}
            &(x_4x_8y_3y_7-x_3x_8y_4y_7-x_4x_7y_3y_8+x_3x_7y_4y_8)f_{15}f_{26}\\
            &+(x_6x_8y_4y_7-x_4x_8y_6y_7-x_6x_7y_4y_8+x_4x_7y_6y_8)f_{15}f_{23}\\
            &+(x_3x_8y_5y_7-x_5x_8y_3y_7-x_3x_7y_5y_8+x_5x_7y_3y_8)f_{14}f_{26}\\
            &+(x_5x_8y_6y_7-x_6x_8y_5y_7-x_5x_7y_6y_8+x_6x_7y_5y_8)f_{14}f_{23}\\
            &+(x_3x_7y_5y_6-x_3x_6y_5y_7-x_3x_5y_6y_7+x_5x_6y_3y_7)f_{12}f_{48}\\
            &+(x_4x_6y_5y_8-x_4x_8y_5y_6-x_5x_6y_4y_8+x_4x_5y_6y_8)f_{12}f_{37}\\
            &+(x_7x_8y_5y_6-x_6x_8y_5y_7-x_5x_7y_6y_8+x_5x_6y_7y_8)f_{12}f_{34} = 0.
        \end{split}
    \end{equation}
    Therefore, one has 
    \begin{equation*}
    \begin{split}
        F &= (x_4x_8y_3y_7-x_3x_8y_4y_7-x_4x_7y_3y_8+x_3x_7y_4y_8)X_{15}X_{26}\\
            &+(x_6x_8y_4y_7-x_4x_8y_6y_7-x_6x_7y_4y_8+x_4x_7y_6y_8)X_{15}X_{23}\\
            &+(x_3x_8y_5y_7-x_5x_8y_3y_7-x_3x_7y_5y_8+x_5x_7y_3y_8)X_{14}X_{26}\\
            &+(x_5x_8y_6y_7-x_6x_8y_5y_7-x_5x_7y_6y_8+x_6x_7y_5y_8)X_{14}X_{23}\\
            &+(x_3x_7y_5y_6-x_3x_6y_5y_7-x_3x_5y_6y_7+x_5x_6y_3y_7)X_{12}X_{48}\\
            &+(x_4x_6y_5y_8-x_4x_8y_5y_6-x_5x_6y_4y_8+x_4x_5y_6y_8)X_{12}X_{37}\\
            &+(x_7x_8y_5y_6-x_6x_8y_5y_7-x_5x_7y_6y_8+x_5x_6y_7y_8)X_{12}X_{34} \in \Lambda_2 \subset \Lambda.
    \end{split}
    \end{equation*}
    From Remark \ref{rem.syz}, it follows that there is no linear relation in the first syzygy of $J_G$. This implies that $F$ is not an element of the module defined by the first syzygy of $J_G$. Therefore, $J_G$ is not a linear type, since $F$ is quadratic homogeneous element.
\end{proof}

\begin{remark}
   In particular, if a unicyclic graph $G$ contains $C_{4,1}$ as an induced subgraph, then $J_G$ is not of linear type.
\end{remark}

Using Macaulay2 \cite{M2}, we computed that if a graph $G$ contains $C_{6,1}$ as an induced subgraph, then $J_G$ is not of $2$ relation type. Thus, one can pose the following question in general:

\begin{question}
    Let $G$ be a graph. If $G$ contains $C_{2n,1}$ as an induced subgraph, then $J_G$ is not of linear type for all $n \geq 2$. In particular, $J_G$ is not of $n-1$ relation type.
\end{question}

\noindent {\bf Acknowledgement.} We thank an anonymous referee for useful comments and suggestions. The first author thanks A. V. Jayanthan for his valuable comments. The first author thank the Scientific and Technological Research Council of Turkey - T\"UB\.{I}TAK (Grant No: 124F113) for the financial support. The second author gratefully acknowledges partial financial support from the Anusandhan National Research Foundation (ANRF), Government of India, via Core Research Grant No. CRG/2023/007668.

\vspace{2mm}

\noindent {\bf Competing Interests.} The authors declare that they have no known competing financial interests or personal relationships that could have appeared to influence the work reported in this paper.

\end{document}